\documentclass[12pt]{amsart}
\usepackage{mathptmx}

\newcommand{\ZZ}{\mathbb{Z}}
\newcommand{\RR}{\mathbb{R}}
\DeclareMathOperator{\Hilb}{Hilb}

\usepackage{fancyvrb}
\RecustomVerbatimEnvironment{Verbatim}{Verbatim}%
{xleftmargin=15pt,baselinestretch=0.9}

\begin{document}
\author{Winfried Bruns and Gesa K\"ampf}

\address{FB Mathematik/Informatik, Universität Osnabrück, 49069
Osnabrück} \email{wbruns@uos.de, gkaempf@uos.de}
\title[A  Macaulay 2 interface for  Normaliz]{A  Macaulay 2 interface for Normaliz}
\maketitle

\vspace{-\baselineskip}
\begin{abstract}
Normaliz is a tool for the computation of Hilbert bases of
normal affine monoids and related tasks. We describe the
Macaulay~2 interface to Normaliz. It makes
Normaliz accessible for operations on monoid algebras in
Macaulay~2 and, conversely, makes Macaulay~2 a
frontend for Normaliz.
\end{abstract}

\section{Introduction}

Normaliz \cite{Nmz} solves a task that can be described from a
geometric and an (equivalent) algebraic perspective. The
objects of discrete geometry on which Normaliz works are
\emph{finitely generated rational cones} $C$, i.e., subsets of
a space $\RR^d$ consisting of the linear combinations
$a_1x_1+\dots+a_nx_n$ of an integral system of generators
$x_1,\dots,x_n\in\ZZ^d$ with nonnegative coefficients $a_i$. If
$C$ is \emph{pointed} $(x,-x\in C \implies x=0$), then the
monoid $M=C\cap\ZZ^d$ has a unique finite minimal system of
generators $\Hilb(M)$, called its \emph{Hilbert basis}. (See
Bruns and Gubeladze \cite{BG} for the terminology and
mathematical background.) The computation of Hilbert bases is
the primary goal of Normaliz. For information on its algorithms
see Bruns and Ichim \cite{BI}.

By the theorem of Minkowski and Weyl, $C$ is a (rational) finitely
generated cone if and only if it is the set of solutions of a linear
homogeneous system of inequalities (with rational coefficients).
Therefore the computation of Hilbert bases is equivalent to solving
such a system of inequalities. For its input, Normaliz accepts
systems of generators as well as systems of inequalities.

Normaliz got its name from the first task for which it was
designed, namely the computation of normalizations of affine
monoids, i.e., finitely generated submonoids
$N=\ZZ_+x_1+\dots+\ZZ_+x_n$ of a lattice $\ZZ^d$. An element
$y\in\ZZ^d$ is \emph{integral} over $N$ if $kx\in N$ for some
$k\in\ZZ$, $k>0$. The \emph{integral closure} $\widehat N$ of
$N$ in $\ZZ^n$ is the monoid of all elements $x\in\ZZ^d$ that
are integral over $N$. Geometrically it is given by $\widehat
N=C\cap \ZZ^d$ for the cone $C=\RR_+x_1+\dots+\RR_+x_n$. The
\emph{normalization} $\bar N$ of $N$ is $\bar N=C\cap\ZZ N$,
and since $\ZZ N$ is isomorphic to a lattice $\ZZ^r$, the
computation of $\bar N$ amounts to a Hilbert basis computation.
Depending on an input parameter, Normaliz computes
$\Hilb(\widehat N)$ or $\Hilb(\bar N)$.

Our terminology has been chosen in such a way that it is
compatible with that of commutative algebra: after the choice
of a field $K$, we can consider the monoid algebra $S=K[N]$,
its integral closure $\widehat S$ in the Laurent polynomial
ring $K[\ZZ^d]$, and its normalization $\bar S$. Since
$\widehat S=K[\widehat N]$ and $\bar S=K[\bar N]$, the
computation of integral closures of monoid algebras is reduced
to the consideration of the underlying monoids.

The package Normaliz provides direct access to Normaliz as well
as ring theoretic functions applied to monomial subalgebras and
monomial ideals in polynomial rings.

\section{Direct access to Normaliz}

The input to Normaliz is composed of matrices of integers. The
interpretation of each matrix depends on a parameter called
\emph{type}.

In types $0$ and $1$ the rows of the input matrix are
interpreted as the generators of a monoid $N$. In type $0$
Normaliz computes $\Hilb(\widehat N)$, and in type $1$ it
computes $\Hilb(\bar N)$ (notation as above).
{\small\begin{Verbatim}
i1 : loadPackage "Normaliz";

i2 : M = matrix {{0,1},{2,3}}:

i3 : normaliz(M,0)

o3 = RationalCone{gen => | 0 1 |                                   }
                         | 1 2 |
                         | 2 3 |
                  inv => HashTable{height 1 elements => 3         }
                                   hilbert basis elements => 3
                                   homogeneous => true
                                   homogeneous weights => (-1, 1)
                                   index => 2
                                   multiplicity => 2
                                   number extreme rays => 2
                                   number support hyperplanes => 2
                                   rank => 2

o3 : RationalCone
\end{Verbatim}
}

The return value of \texttt{normaliz} is an object of type
\texttt{rationalCone} defined by the package. Such an object is
a hash table whose components are matrices and a hash table.
The matrices are named after the suffixes of the output files
of Normaliz from which the corresponding matrix is read.
This example represents the minimal content of the
\texttt{rationalCone} returned, namely the Hilbert basis in
\texttt{gen} and the numerical invariants in \texttt{inv}.

The amount of data returned can be increased by the option
\texttt{allComputations}:
{\small\begin{Verbatim}
i4 : normaliz(allComputations => true,M,0)

o4 = RationalCone{cgr => 0                                         }
                  equ => 0
                  gen => | 0 1 |
                         | 1 2 |
                         | 2 3 |
                  inv => HashTable{height 1 elements => 3         }
                                   hilbert basis elements => 3
                                   homogeneous => true
                                   homogeneous weights => (-1, 1)
                                   index => 2
                                   multiplicity => 2
                                   number extreme rays => 2
                                   number support hyperplanes => 2
                                   rank => 2
                  sup => | -3 2 |
                         | 1  0 |
                  typ => | 2 0 |
                         | 1 1 |
                         | 0 2 |

o4 : RationalCone
\end{Verbatim}
}
In addition to the data above, the rational cone will now contain
four more components: \texttt{cgr}, \texttt{equ}, \texttt{sup},and
\texttt{typ}. The matrix \verb+arc#"sup"+ contains the linear forms
defining the cone. Each row $(\alpha_1,\dots,\allowbreak \alpha_d)$
represents an inequality $\alpha_1\xi_1+\dots+\alpha_d\xi_d\ge 0$
for the vectors $(\xi_1,\dots,\xi_d)$ of the cone. The matrix
\verb+arc#"typ"+ contains the values of these linear forms on the
Hilbert basis elements. The matrices \verb+arc#"cgr"+ and
\verb+arc#"equ"+ are empty in this example: they contain the
equations and congruences that together with the inequalities define
the cone and the lattice, respectively, and in this example the
inequalities are sufficient.

The output data of Normaliz are further increased by
\texttt{setNmzOption("allf",true)}. For the complete list of
data and their interpretation see the Normaliz documentation or
the package online help.

Normaliz has 6 more input types. Types 4, 5 and 6 allow the
user to specify the cone and the lattice by homogeneous
diophantine equations, inequalities and congruenences. A type 4
matrix is interpreted as a system of inequalities, defining the
cone $C$ as just explained in connection with the matrix
\texttt{"sup"}. The rows $(\alpha_1,\dots,\alpha_d)$ of a type
5 matrix are considered as equations
$\alpha_1\xi_1+\dots+\alpha_d\xi_d=0$. In type 6 the matrix is
interpreted as a system of homogeneous congruences, and $\ZZ^d$
is then replaced by the lattice of solutions. Types 4, 5 and 6
can be combined, and therefore the function \verb+normaliz+ has
a variant in which the argument of the function is a list
consisting of pairs \verb+(matrix, type)+. (The default type 4
matrix is the unit matrix, defining the positive orthant.)

For these types, \texttt{setNMzOption("dual",true)} chooses an
alternative algorithm. It is often better than the
triangulation based standard algorithm of Normaliz.

The input types 2 and 3 are variants of type 0. In type 2 the
rows $x_1,\dots,x_n$ of the matrix are interpreted as the
vertices of a lattice polytope, and type $0$ is applied to the
cone generated by $x_i'=(x_i,1)\in \RR^{d+1}$, $i=1,\dots,n$.
Types 3 and 10 have a ring theoretic flavor. They will be
explained in the next section. (Types 7, 8 and 9 are reserved for
future extensions.)

Via suitable options one can restrict the data that Normaliz
computes (see online help or Normaliz documentation). The
extension \texttt{setNmzOption("hilb",true)} is more important.
It asks Normaliz to find the Hilbert series and polynomial of
the (algebra over the) integral closure or normalization
computed, provided a homogeneity condition is satisfied: there
is an integral linear form $\lambda$ on $\ZZ^d$ or $\ZZ N$,
respectively, such that $\lambda(x)=1$ for the extreme integral
generators $x$ of the cone. For the introductory example, the
data in \verb+arc#"inv"+ now contain
{\small
\begin{Verbatim}
i5 : setNmzOption("hilb",true);

i6 : normaliz(M,0)

o6 = RationalCone{gen => | 0 1 |                                   }
                         | 1 2 |
                         | 2 3 |
                  inv => HashTable{h-vector => (1, 1)             }
                                   height 1 elements => 3
                                   hilbert basis elements => 3
                                   hilbert polynomial => (1, 2)
                                   homogeneous => true
                                   homogeneous weights => (-1, 1)
                                   index => 2
                                   multiplicity => 2
                                   number extreme rays => 2
                                   number support hyperplanes => 2
                                   rank => 2

o6 : RationalCone
\end{Verbatim}
}
The $h$-vector represents the numerator polynomial $1+t$ of the Hilbert
series, and the Hilbert polynomial is $1+2k$ where $k$ denotes the degree.

\section{Ring theoretic functions}

For the ring theoretic functions the package introduces the class
\texttt{monomialSubalgebra}. It is a subclass of \texttt{Ring}. A
monomial subalgebra is created as follows:
{\small\begin{Verbatim}
i7 : R=ZZ/17[x,y];

i8 : S = createMonomialSubalgebra {x,x^2*y^3}

     ZZ     2 3
o8 = --[x, x y ]
     17

o8 : monomial subalgebra of R
\end{Verbatim}
}
This creates the monomial $K$-subalgebra $S$ generated by the
monomials $x,x^2y^3$ in the polynomial ring $R=K[x,y]$
over the base field $K$.

The  functions \verb+normalToricRing+ and \verb+intclToricRing+
have monomial subalgebras as input (or just lists of
monomials): {\small
\begin{Verbatim}
i9 : intclToricRing S

      ZZ          2 3
o9 = --[x, x*y, x y ]
      17

o9 : monomial subalgebra of R
\end{Verbatim}
}

\texttt{normalToricRing} returns the normalization of $S$ as a
monomial subalgebra of $R$ ($R$ always contains the normalization),
whereas \texttt{intclToricRing} returns the integral closure of $S$
in $R$ (or, equivalently, in the field of fractions of $R$). The
cache of the returned monomial subalgebras contains the rational cone
computed by Normaliz.

The function \texttt{intclMonIdeal} has a monomial ideal $I$ as its
input. It computes the integral closure $\bar I$ of $I$ and the
normalization of the Rees algebra of $I$. Consequently the return
value is a sequence containing $\bar I$ (of type \texttt{ideal}) and
a monomial subalgebra. Note that the Rees algebra $R[It]$ and its
normalization live in the extended polynomial ring $R[t]$. The
function creates this extended polynomial ring, choosing an
available name for the auxiliary indeterminate $t$. If $R$ itself is
of type $R'[t]$, $I$ is a monomial ideal in $R'$ and the
normalization of the Rees algebra $R'[It]$ is to be computed , then
the user can indicate this fact by adding the name of $t$ to the input.

The function \texttt{intclToricRing} calls Normaliz in type 3,
made exactly for the computation of normalizations of Rees
algebras.

The function \texttt{normalToricRing} has a variant in which the
input parameter is an ideal $I$ consisting of binomaials $X^a-X^b$.
It has a unique minimal prime ideal $P$ generated by binomials of
the same type, and the function returns the normalization of $R/P$
embedded into a newly created polynomial ring of the same Krull
dimension. (In general there is no canonical choice for such an
embedding.)
{\small\begin{Verbatim}
i10 : R = ZZ/37[x,y,z,w];

i11 : I = ideal(x*w-y*z,x*z-y^2);

o11 : Ideal of R

i12 : normalToricRing(I,t)

      ZZ  3     2   2     3
o12 = --[t , t t , t t , t ]
      37  2   1 2   1 2   1

                             ZZ
o12 : monomial subalgebra of --[t , t ]
                             37  1   2
\end{Verbatim}
}
This function uses type 10 of \texttt{normaliz}, created exactly for
this purpose.

There are further ring theoretic functions in the library:
\texttt{intersectionValRings},
\texttt{intersectionValRingIdeals}, \texttt{torusInvariants},
\texttt{finiteDiagInvariants} and\linebreak[4]
\texttt{diagInvariants}. The first two compute intersections of
monomial valuation rings and ideals with the polynomial ring
$R$ whereas the last three compute the rings of invariants of a
diagonal torus action on $R$, a diagonal finite group action,
or a diagonal group action in general.

\section{Miscellanea}

Macaulay 2 and Normaliz exchange data via hard disk files. By
default the package handles the files behind the scenes, and the
user need not care about them. However, the user can take over
command of the file handling by specifying a file name and a path to
the directory where the files are to be stored. The package provides
functions for writing and reading Normaliz files directly. See the
online help for details.

The standard integer precision of Normaliz is 64 bit
(corresponding to the C integer type \texttt{long long}).
Already in small dimensions this may not be sufficient. In that
case one can choose the indefinite precision executable by
setting \texttt{nmzVersion="normbig"}. This choice typically
increases the computation time by a factor of 5. It is less
time consuming to use \texttt{setNmzOption("errorcheck",true)}
in order to control the arithmetic of the 64 bit computation.
For further functions and options we refer the reader to the
online help.
\medskip

\emph{Acknowledgement} The authors are grateful to Dan Grayson
for his support in the  development of the Normaliz package.

\end{document}